\documentclass[final,twoside,11pt,final]{entics}
\usepackage{enticsmacro}
\usepackage{url}
\usepackage{xspace}
\usepackage{tikz-cd}
\usepackage{amsmath}
\usepackage{mathpartir}
\usepackage{scalerel}
\usepackage{subcaption}
\newcommand{\catt}{\texttt{CaTT}\xspace}
\newcommand{\gsett}{\texttt{GSeTT}\xspace}

\newcommand{\N}{\ensuremath{\mathbb{N}}}
\newcommand{\T}{\ensuremath{\mathbb{T}}}
\newcommand{\disk}{\ensuremath{\mathbb{D}}}
\newcommand{\sphere}{\ensuremath{\mathbb{S}}}
\newcommand{\terminal}{\ensuremath{\mathbf{1}}}

\newcommand{\syn}{\ensuremath{\mathcal{S}}}
\DeclareMathOperator*{\bigHash}{\scalerel*{\#}{\sum}}

\DeclareMathOperator{\Var}{Var}
\newcommand{\vars}{\ensuremath{\mathcal{V}}}
\newcommand{\sub}[1]{\ensuremath{\langle #1 \rangle}}
\newcommand{\obj}{\ast}
\newcommand{\arr}[3][]{#2 \to_{#1} #3}
\newcommand{\emptycontext}{\ensuremath{\emptyset}}
\newcommand{\vdashps}{\ensuremath{\vdash_{\textsf{ps}}}}

\newcommand{\cat}{\mathcal}

\DeclareMathOperator{\PsCtx}{PsCtx}

\DeclareMathOperator{\Bat}{Bat}
\DeclareMathOperator{\Set}{Set}
\DeclareMathOperator{\Comp}{Comp}
\DeclareMathOperator{\gSet}{gSet}

\DeclareMathOperator{\Fam}{Fam}

\DeclareMathOperator{\List}{List}
\DeclareMathOperator{\Sol}{Sol}
\DeclareMathOperator{\Zig}{Zig}
\DeclareMathOperator{\Card}{Card}
\DeclareMathOperator{\Pos}{Pos}
\DeclareMathOperator{\Free}{Free}
\DeclareMathOperator{\Cell}{Cell}
\DeclareMathOperator{\Sphere}{Sphere}
\DeclareMathOperator{\Term}{Tm}
\DeclareMathOperator{\Type}{Ty}
\DeclareMathOperator{\tr}{tr}
\DeclareMathOperator{\Tr}{Tr}
\DeclareMathOperator{\op}{op}
\DeclareMathOperator{\Full}{Full}
\DeclareMathOperator{\var}{var}
\DeclareMathOperator{\coh}{coh}
\DeclareMathOperator{\supp}{supp}
\DeclareMathOperator{\rmlast}{rm}
\DeclareMathOperator{\image}{im}
\DeclareMathOperator{\Tree}{Tree}

\DeclareMathOperator{\src}{src}
\DeclareMathOperator{\tgt}{tgt}
\DeclareMathOperator{\br}{br}
\DeclareMathOperator{\bdry}{bdry}
\DeclareMathOperator{\pr}{pr}
\DeclareMathOperator{\id}{id}

\def\conf{MFPS 2024}
\def\myconf{mfps40} 	
\volume{4}			
\def\volu{4}			
\def\papno{5}			
\def\mydoi{14675}%

\def\lastname{Benjamin, Markakis, Sarti}
\def\runtitle{CaTT contexts are finite computads}
\def\copynames{T. Benjamin, I. Markakis, C. Sarti}

\usepackage{mathpartir}

\makeatletter

\NewDocumentCommand \mpr@label { m }{%
  \def\@currentlabelname{#1}%
   \phantomsection 
}

\DeclareDocumentCommand \inferdef {m m m}{%
   \inferrule*[vcenter,right=#1]{#2}{#3}
   \mpr@label{\textsc{#1}}
}%

\makeatother

\newcommand{\ruleref}[1]{Rule \nameref{#1}}

\begin{document}
\begin{frontmatter}
  \title{CaTT Contexts are Finite Computads\thanksref{ALL}}
  \thanks[ALL]{We would like to thank Jamie Vicary for his feedback.}
  \author{Thibaut Benjamin\thanksref{a}\thanksref{tbemail}}
  \author{Ioannis Markakis\thanksref{a}\thanksref{imemail}\thanksref{b}}
  \author{Chiara Sarti\thanksref{a}\thanksref{csemail}}
  \address[a]{Department of Computer and Technology\\ University of Cambridge\\
    Cambridge, United Kingdom}
  \thanks[tbemail]{Email: \href{mailto:tjb201@cam.ac.uk}
    {\texttt{\normalshape tjb201@cam.ac.uk}}}
  \thanks[imemail]{Email:  \href{mailto:ioannis.markakis@cl.cam.ac.uk}
    {\texttt{\normalshape ioannis.markakis@cl.cam.ac.uk}}}
  \thanks[csemail]{Email:  \href{mailto:cs2197@cam.ac.uk}
    {\texttt{\normalshape cs2197@cam.ac.uk}}}
  \thanks[b]{The author would like to acknowledge funding by the Onassis
    Foundation - Scholarship ID: {F ZQ 039-1/2020-2021}}

\begin{abstract}
  Two novel descriptions of weak \(\omega\)-categories have been recently
  proposed, using type-theoretic ideas. The first one is the dependent type
  theory \catt whose models are \(\omega\)-categories. The second is a recursive
  description of a category of computads together with an adjunction to globular
  sets, such that the algebras for the induced monad are again
  \(\omega\)-categories. We compare the two descriptions by showing that there
  exits a fully faithful morphism of categories with families from the syntactic
  category of \catt to the opposite of the category of computads, which gives an
  equivalence on the subcategory of finite computads. We derive a more direct
  connection between the category of models of \catt and the category of
  algebras for the monad on globular sets, induced by the adjunction with
  computads.
\end{abstract}
\end{frontmatter}

\section{Introduction}\label{sec:introduction}

Higher categories have recently found applications in theoretical computer
science and other fields. They describe the iterated structure of identity types
in Martin-L\"{o}f type
theory~\cite{altenkirch_syntactical_2012,lumsdaine_weak_2009,vandenberg_types_2011},
a fact that has been instrumental in the emergence of homotopy type
theory~\cite{theunivalentfoundationsprogram_homotopy_2013}. They are also
crucial in higher-dimensional rewriting~\cite{mimram_3dimensional_2014},
homotopy theory ~\cite{garner_homomorphisms_2010,henry_homotopy_2023} and
homology theory~\cite{lafont_polygraphic_2009}. There exists a wide range of
flavours of higher categories, differentiated by the shapes of their higher
morphisms and the ways that they are composed. For a review of the various
definitions, we refer to the work of Leinster~\cite{leinster_survey_2001}.

In this paper, we focus on globular weak \(\omega\)-categories, a particular
model of higher categories. Different definitions of weak \(\omega\)-categories
have been given by Batanin and
Leinster~\cite{batanin_monoidal_1998,leinster_higher_2004}, and by Grothendieck
and
Maltsiniotis~\cite{grothendieck_pursuing_1983,maltsiniotis_grothendieck_2010}.
More recently, two syntactic descriptions of weak
\(\omega\)-categories have been proposed. One takes the form of a
dependent type theory called \catt, introduced by Finster and
Mimram~\cite{finster_typetheoretical_2017}. The other one gives a description of
the free \(\omega\)\-category monad using inductive types, which was introduced
by Dean et al.~\cite{dean_computads_2022} and is given via an adjunction to a
category of \emph{computads}. The equivalence between models of the type theory
\catt and algebras of the free \(\omega\)\-category monad of Dean et al. has
been established by a chain of results. Benjamin, Finster and
Mimram~\cite{benjamin_globular_2021} show that models of \catt are equivalent to
Grothendieck-Maltsiniotis weak \(\omega\)-categories. Dean et
al.~\cite{dean_computads_2022} show that their monad is the one of Leinster.
Finally, the equivalence of the theories of Batanin-Leinster and
Grothendieck-Maltsiniotis has been established by Ara~\cite{ara_infty_2010} and
Bourke~\cite{bourke_iterated_2020}. In this paper, we establish a more direct
comparison between the two syntactic descriptions by showing that the contexts
of \catt correspond to finite computads.

The notion of \(\omega\)\-categories is an extension of the notion of globular
sets. The data of a globular sets is a set of cells of various dimensions,
equipped with source and target operations. Those operations allows us to
visualise cells as directed analogues of disks of various dimensions, as
depicted below in low dimensions:
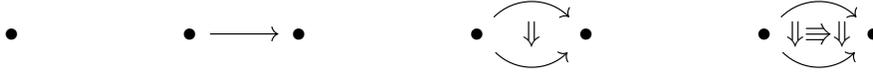
\[
  \begin{tikzcd}
    \bullet
    & & \bullet\ar[r] & \bullet
    & & \bullet\ar[r, bend left = 45]\ar[r, bend right= 45]\ar[r, phantom, "\Downarrow"]
    & \bullet
    & & \bullet\ar[r, bend left= 45]\ar[r, bend right= 45]
      \ar[r, phantom, "\Downarrow\Rrightarrow\Downarrow"]
    & \bullet
  \end{tikzcd}
\]
An \(\omega\)-category is a globular set \(G\) equipped with ways to compose
certain configurations of cells. Those operations generalise for instance the
composition operations of ordinary categories and bicategories, in the sense
that the following diagrams of cells admit a composite:
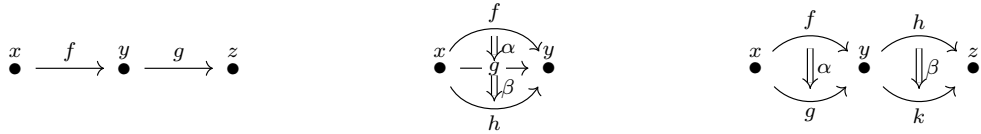
\begin{align*}
  \begin{tikzcd}[ampersand replacement=\&]
    \overset{x}{\bullet}\ar[r,"f"]
    \& \overset{y}{\bullet}\ar[r,"g"]
    \& \overset{z}{\bullet}
  \end{tikzcd} &&
  \begin{tikzcd}[ampersand replacement=\&]
    \vphantom{\bullet}\smash{\overset{x}{\bullet}} \&
    \vphantom{\bullet}\smash{\overset{y}{\bullet}}
    \arrow[""{name=0, anchor=center, inner sep=0}, "f", bend left = 60, from=1-1, to=1-2]
    \arrow[""{name=1, anchor=center, inner sep=0}, "h"', bend right = 60, from=1-1, to=1-2]
    \arrow[""{name=2, anchor=center, inner sep=0}, "g"{description}, from=1-1, to=1-2]
    \arrow["\alpha", shorten <=2pt, shorten >=2pt, Rightarrow, from=0, to=2]
    \arrow["\beta", shorten <=2pt, shorten >=2pt, Rightarrow, from=2, to=1]
  \end{tikzcd} &&
  \begin{tikzcd}[ampersand replacement=\&]
    \vphantom{\bullet}\smash{\overset{x}{\bullet}}
    \ar[r,"f", bend left = 45,""{name=0, anchor=center, inner sep=0}]
    \ar[r,"g"', bend right= 45,""{name=1, anchor=center, inner sep=0}]
    \& \vphantom{\bullet}\smash{\overset{y}{\bullet}}
    \ar[r,"h", bend left= 45,""{name=2, anchor=center, inner sep=0}]
    \ar[r,"k"', bend right= 45,""{name=3, anchor=center, inner sep=0}]
    \& \vphantom{\bullet}\smash{\overset{z}{\bullet}}
    \ar["\alpha", shorten <=4pt, shorten >=4pt, Rightarrow, from=0, to=1]
    \ar["\beta", shorten <=4pt, shorten >=4pt, Rightarrow, from=2, to=3]
  \end{tikzcd}
\end{align*}
One may expect from ordinary category theory that those operations are strictly
associative and unital. This leads to the notions of strict
\(\omega\)\-categories. However, for an arbitrary \(\omega\)\-category,
associativity and unitality only hold in a weaker sense. More precisely, there
exist operations producing higher-dimensional \emph{coherence} cells witnessing
those axioms. For example, there is an operation, depicted below, producing a
witness cell for the associativity of the composition of three \(1\)-cells:
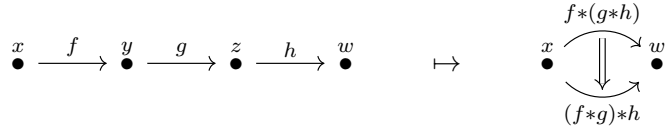
\[
  \begin{tikzcd}
    \overset{x}{\bullet}\ar[r,"f"]
    & \overset{y}{\bullet}\ar[r,"g"]
    & \overset{z}{\bullet}\ar[r,"h"]
    & \overset{w}{\bullet}
  \end{tikzcd}
  \qquad \mapsto \qquad
  \begin{tikzcd}
    \vphantom{\bullet}\smash{\overset{x}{\bullet}}
    \ar[r, bend left = 45,"f*(g*h)", ""{name=0, anchor=center, inner sep=0}]
    \ar[r, bend right = 45,"(f*g)*h"', ""{name=1, anchor=center, inner sep=0}]
    \ar[from = 0, to = 1, shorten <=2pt, shorten >=2pt, Rightarrow ]
    & \vphantom{\bullet}\smash{\overset{w}{\bullet}}
  \end{tikzcd}
\]
This leads to an infinite tower of operations, where higher witness cells relate
composites of lower-dimensional witness cells. For example, there are cells
relating the two sides of the pentagon and triangle equations of monoidal
categories.

The dependent type theory \catt gives a syntactic way to systematically track
those operations. Types in \catt combinatorially encode directed spheres, while
terms encode directed disks. The terms are produced by a unique term constructor
\(\coh\) parametrised by a context \(\Gamma\) of a particular kind called a
\emph{pasting context} together with a type \(\Gamma \vdash A\) satisfying some
side conditions. Pasting contexts encode the arities of the composition and
coherence operations of \(\omega\)\-categories, i.e. the configurations of
cells that may be composed. When the side condition is satisfied, the
term constructor \(\coh_{\Gamma,A}\) produces a canonical inhabitant for the
type \(A\). This inhabitant can be seen as the composite cell of the variables
of \(\Gamma\), or a coherence cell relating the source and target of \(A\),
depending on the dimension of \(A\). We note that this type theory has no
equality between terms, due to the lack of equations between the operations of
the \(\omega\)\-categories it models.

Contexts in the type theory \catt are finite lists of generators for
\(\omega\)\-categories. An alternative way to freely generate
\(\omega\)\-categories, due to Street~\cite{street_limits_1976},
Burroni~\cite{burroni_higherdimensional_1993} and
Batanin~\cite{batanin_computads_1998}, is that of \emph{computads} or
\emph{polygraphs}. Computads are sets of generators stratified by dimension,
together with attaching functions assigning to each cell a source and target,
which are lower-dimensional cells of the free \(\omega\)\-category generated by
the computad. More precisely, computads and the free \(\omega\)\-category on
them are defined by mutual induction on their dimension. Inspired by the type
theory \catt, Dean et al.~\cite{dean_computads_2022} recently gave a new
presentation of the free \(\omega\)\-category on a computad using inductive
types. In this presentation, the sets of cells of the free \(\omega\)\-category
are inductively generated by a constructor \(\coh\) with similar arguments as
the term constructor of \catt.

Our main contribution is showing that the syntactic category of \catt is
equivalent to the opposite of the category of finite
computads. We do this by first equipping the opposite of the
category of computads with a structure of category with families. We then
establish a bijection between pasting contexts and \emph{Batanin trees}, the
arities of the constructor \(\coh\) in the case of computads. We construct
recursively a morphism of categories with families
\(R \colon \syn_{\catt} \to \Comp^{\op}\) from the syntactic category of \catt
to the opposite of the category of computads, comparing types \(A\) satisfying
the side condition of the term constructor to the \emph{full spheres} of the
cell constructor. Finally, in our main Theorem~\ref{thm:main}, we shown that
this morphism is fully faithful and essentially surjective on the subcategory of
finite computads, establishing the claimed equivalence between contexts of \catt
and finite computads.

\paragraph{Overview of the paper.}
In Section~\ref{sec:catt}, we give an overview of categories with families and
the dependent type theories \gsett and \catt. In Section~\ref{sec:computads}, we
present the definition of globular sets and computads, and the adjunction
between them. Then in Section~\ref{sec:globular-gsett}, we briefly describe the
equivalence between the syntactic category of \gsett and finite globular sets,
established by Benjamin, Finster, Mimram~\cite{benjamin_globular_2021}. In
Section~\ref{sec:pasting-batanin}, we compare the pasting contexts of \catt to
the Batanin trees. Section~\ref{sec:main} contains our main result, the
comparison of the syntactic category of \catt and the category of finite
computads. In Section~\ref{sec:models}, we discuss briefly how our result lets
us obtain a more direct proof of the equivalence of models of \catt and algebras
of the monad on globular sets induced by the adjunction with computads.

\section{The dependent type theory CaTT} \label{sec:catt}
\subsection{Dependent type theories and categories with families}
\begin{figure}[b]
  \centering
  \fbox{
    \begin{minipage}{.95\linewidth}
      \begin{mathpar}
        \inferdef{(ec)}{\null}{\emptycontext\vdash}\label{rule:a} \and
        \inferdef{(var)}{(x:A)\in\Gamma}{\Gamma\vdash \var x:A}
        \and \inferdef{(es)}{\Gamma\vdash}{\Gamma\vdash \sub{}:\emptycontext}\\
        \inferdef{(cc)}{\Gamma\vdash A \\ x\notin \Var(\Gamma)}{(\Gamma,x:A)
          \vdash} \and \inferdef{(sc)}
        {\Delta\vdash \gamma:\Gamma \\
          (\Gamma, x:A)\vdash \\ \Delta\vdash t:A[\gamma]} {\Delta
          \vdash\sub{\gamma,x\mapsto t} : (\Gamma,x : A)}
      \end{mathpar}
    \end{minipage}
  }
  \caption{Context and substitution rules}
  \label{fig:rules-basic}
\end{figure}

The dependent type theory \(\catt\), introduced by Finster and
Mimram~\cite{finster_typetheoretical_2017}, is a syntactic way to work with weak
\(\omega\)-categories. To introduce a dependent type theory, we consider a fixed
set of variables \(\vars\) that we assume countably infinite, and define
mutually inductively the following syntactic objects
\begin{itemize}
\item Contexts are lists of pairs \((x:A)\), where \(x\in\vars\) and \(A\) is a
  type (empty list denoted \(\emptycontext\))
\item Substitutions are lists of \(\sub{x\mapsto t}\) where \(x\in\vars\) and
  \(t\) is a term (empty list denoted \(\sub{}\))
\item Types are expressions built from type constructors
\item Terms are expressions either of the form \(\var x\) where \(x\) is a
  variable, or built from term constructors
\end{itemize}
We note that only some of those syntactic objects will be considered valid in
the type theory, as explained below.

A substitution \(\gamma\) may act on a type \(A\) (resp.~on a term \(t\)),
producing a type \(A[\gamma]\) (resp.~a term \(t[\gamma]\)). The definition may
vary depending on the type theory, but we require that for every variable \(x\),
the term \((\var x)[\gamma]\) is the last term associated to \(x\) in \(\gamma\)
if it exists, and \(\var x\) otherwise. A set of variables is associated to
contexts, terms, types and substitutions, denoted \(\Var\), and we always define
\begin{align*}
  \Var(\emptycontext) &= \emptyset
  & \Var(\Gamma,x:A) &= \Var(\Gamma)\cup\{x\} \\
  \Var(\sub{}) &= \emptyset
  & \Var(\sub{\gamma,x\mapsto t}) &= \Var(\gamma)\cup\Var(t) \\
  \Var(\var x) &= \{x\}
\end{align*}
Additionally, contexts, substitutions, types and
terms are subjects to well-formedness conditions that are expressed in the form
of the following judgements, presented here along with their intended meaning.
\[
  \begin{array}{r@{\qquad}l}
  \Gamma\vdash & \text{\(\Gamma\) is a valid context}\\
  \Delta\vdash \gamma:\Gamma & \text{\(\gamma\) is a valid substitution from
                               \(\Delta\) to \(\Gamma\)}
\end{array}
\qquad\qquad
\begin{array}{r@{\qquad}l}
  \Gamma\vdash A & \text{\(A\) is a valid type in \(\Gamma\)} \\
  \Gamma\vdash t:A & \text{\(t\) is a type of type \(A\) in \(\Gamma\)}
\end{array}
\]
We also require all dependent type theories to satisfy the context and
substitutions forming rules, as well as the variable forming rules presented in
Fig.~\ref{fig:rules-basic}.
Additionally, in the theories that we consider, the action of substitution on
terms defines composition of substitutions which is associative and unital.

This lets us consider the \emph{syntactic category}
\(\syn_{\T}\) associated to a dependent type theory \(\T\), whose objects are
the valid contexts and whose morphisms \(\Delta\to \Gamma\) are the
substitutions \(\Delta\vdash \gamma:\Gamma\). the
syntactic category carries a structure of category with families, in the sense of
Dybjer~\cite{dybjer_internal_1996}. We recall that the category \(\Fam\) of
\emph{families of sets} is the category with objects set-indexed families of
sets \((A_i)_{i\in I}\) and morphisms are pairs
\((\phi,f_i)\colon (A_i)_{i\in I}\to (B_j)_{j\in J}\) of a function
\( \phi \colon I\to J\) together with a family of functions
\(f_i\colon A_i\to B_{\phi(i)}\) for every \(i\in I\). A \emph{category with
families} \(\cat C\) is a category with a chosen terminal object
\(\terminal\), equipped with
a \(\Fam\)\-valued presheaf \(T \colon {\cat C}^{\op}\to \Fam\), denoted by
\[
  T(\Gamma) = (\Term(\Gamma;A))_{A\in \Type\Gamma}.
\]
Additionally, there exists an extension operation associating to every pair of
an object \(\Gamma \in \cat C\) and \(A\in \Type \Gamma\), an object
\((\Gamma.A)\) together with a projection
\(\pi_{\Gamma,A}:(\Gamma.A) \to \Gamma\) and an element
\(v_{\Gamma,A} \in \Term(\Gamma;\Type(\pi_{\Gamma,A})A)\) satisfying the
following universal property: For every object \(\Delta\) with a map
\(\gamma:\Delta\to \Gamma\) and an element \(u\in\Term(\Delta;\Type(\gamma)A)\),
there exists a unique map \(\sub{\gamma,u}:\Delta\to(\Gamma.A)\) such that
\begin{align*}
  \pi_{\Gamma,A}\sub{\gamma,u} &= \gamma &
  \Term(\sub{\gamma,u})(v_{\Gamma,A}) &= u.
\end{align*}
A morphism of categories with families \(F \colon \cat C\to \cat D\) consists of
a functor \(\cat C\to \cat D\) preserving the chosen terminal object, together
with a natural transformation
\(F^T \colon T^\cat{C}\Rightarrow T^{\cat D} \circ F\) such that the extension
operation is preserved on the nose. By definition of the category of families,
we see that \(F^T\) amounts to functions
\begin{align*}
  F^{\Type}_\Gamma &\colon \Type^{\cat C}(\Gamma) \to \Type^{\cat D}(F\Gamma) &
  F^{\Term}_{\Gamma,A} &\colon
    \Term^{\cat C}(\Gamma;A) \to \Type^{\cat D}(F\Gamma; F^{\Type}_\Gamma A)
\end{align*}
satisfying the expected naturality conditions.

The structure of category with families for the syntactic category of a
dependent type theory is given by choosing \(\Type(\Gamma)\) to be the valid
types in context \(\Gamma\), choosing \(\Term(\Gamma;A)\) to be the set of terms
of type \(A\) in context \(\Gamma\), and the extension operation to be the
context extension of the theory. The universal projection
\(\pi_{\Gamma,A}:(\Gamma.A)\to \Gamma\) is the
obvious weakening and the universal element \(v_{\Gamma,A}\) is the last
variable of the context. We have formalised~\cite{benjamin_formalisation_2022}
the category with families structure of the syntactic category for the two
dependent type theories \gsett and
\catt that we consider here. For the sake of readability, we introduce the
dependent type theories using named variables, however we consider the contexts
in the syntactic category to be quotiented by \(\alpha\)-renaming, which can be
achieved in practice by presenting the type theory with De Bruijn indices for
instance~\cite{debruijn_lambda_1972}. We make sure to present the theories in a
way that makes translating to De Bruijn indices trivial.

\subsection{The dependent type theories GSeTT and CaTT}
\begin{figure}
  \centering
  \fbox{
    \begin{minipage}{.95\linewidth}
      \begin{mathpar}
        \obj[\gamma] = \obj
        \and (\arr[A]{t}{u})[\gamma] = \arr[{A[\gamma]}]{t[\gamma]}{u[\gamma]}\\
        \Var{(\obj)} = \emptyset \and \Var(\arr[A]{t}{u}) = \Var(A)\cup
        \Var(t)\cup\Var(u)\\
        \inferdef{(obj)}{\Gamma\vdash}{\Gamma\vdash \obj} \and
        \inferdef{(arr)}{\Gamma\vdash A \\ \Gamma\vdash t:A \\ \Gamma\vdash
          u:A}{\Gamma\vdash \arr[A]{t}{u}}
      \end{mathpar}
    \end{minipage}
  }
  \caption{Definition of the type theory \gsett}
  \label{fig:rules-gsett}
\end{figure}
\begin{figure}
  \centering
  \fbox{
    \begin{minipage}{.95\linewidth}
      \begin{mathpar}
        \inferdef{(pss)}{\null}{(x:\obj)\vdashps x:\obj}\label{rule:pss} \and
        \inferdef{(psd)}{\Gamma\vdashps f:\arr[A]{x}{y}}{\Gamma\vdashps y:A}
        \label{rule:psd}\\
        \inferdef{(pse)}{\Gamma\vdashps x:A \\ y,f\notin \Var(\Gamma)}
        {\Gamma,y:A,f:\arr[A]{x}{y}\vdashps f:\arr[A]{x}{y}}\label{rule:pse}
        \and \inferdef{(ps)}{\Gamma\vdashps
          x:\obj}{\Gamma\vdashps}\label{rule:ps}
      \end{mathpar}
    \end{minipage}
  }
  \caption{Derivation of pasting schemes}
  \label{fig:rules-ps}
\end{figure}
\begin{figure}
  \centering
  \fbox{
    \begin{minipage}{.95\linewidth}
      \begin{mathpar}
        \coh_{\Gamma,A}[\gamma][\delta] = \coh_{\Gamma,A}[\gamma\circ\delta]
        \and
        \Var(\coh_{\Gamma,A}[\gamma]) = \Var(\gamma)\\
        \inferdef{(coh)}{\Gamma\vdashps \\ \Gamma\vdash u:A \\ \Gamma\vdash v:A \\
          \Delta\vdash \gamma:\Gamma}{\Delta\vdash
          \coh_{\Gamma,\arr[A]{u}{v}}[\gamma]:(\arr[A]{u}{v})[\gamma]}
        \label{rule:coh} \\
        \text{Where \ruleref{rule:coh} is subject to one of the following
          side conditions} \\
        \dim \Gamma \ge 1 \quad\text{and}\quad \begin{cases}
          \Var (s^\Gamma_{\dim \Gamma -1}) = \Var u \cup \Var A\\
          \Var (t^\Gamma_{\dim \Gamma -1}) = \Var v \cup \Var A
        \end{cases}
        \and \text{or} \and
        \begin{cases}
          \Var(\Gamma) = \Var u \cup \Var A\\
          \Var(\Gamma) = \Var v \cup \Var A
        \end{cases}
      \end{mathpar}
    \end{minipage}
  }
  \caption{Definition of the type theory \catt}
  \label{fig:rules-catt}
\end{figure}

We first introduce the dependent type theory \gsett, a dependent type theory
describing globular sets. This theory has no term constructors, and two type
constructors \(\obj\), taking no argument, and \(\arr[]{}{}\), taking as a
arguments a type \(A\) and two terms \(u\) and \(v\) and producing the type
\(\arr[A]{t}{u}\). This type theory is completely determined by the definitions
and
rules given in Fig.~\ref{fig:rules-gsett}. The dependent type theory \catt
is an extension of \gsett, with one term constructor \(\coh\) taking as
arguments a context \(\Gamma\), a type \(A\) and a substitution \(\gamma\) and
producing the term \(\coh_{\Gamma,A}[\gamma]\). We note that the brackets here
are part of the syntax, and do not denote the application of a substitution.
The introduction rule for this
term constructor requires the use of the two following auxiliary judgements,
subject to the rules presented in Fig.~\ref{fig:rules-ps}
\[
  \begin{array}{r@{\quad}l@{\qquad\qquad}r@{\quad}l}
    \Gamma\vdashps & \text{\(\Gamma\) is a pasting scheme}
    &  \Gamma\vdashps x:A & \text{auxiliary judgement}
  \end{array}
\]
Moreover, we need the following auxiliary definitions. The \emph{dimension} of a
type and a context are defined recursively by
\begin{align*}
  \dim(\obj)&= -1 & \dim (\arr[A]{t}{u}) &= \dim{A} + 1 \\
  \dim(\emptycontext) &= -1 & \dim(\Gamma,x:A) &= \max(\dim(\Gamma), \dim A +1)
\end{align*}
For every pasting context \(\Gamma\vdashps\) and every \(k\in\N\), we
define two contexts \(\partial^\pm_k(\Gamma)\) of \gsett by induction on the
derivation, where \(\rmlast\) is the function removing the last element of a
list and \(\epsilon\in \{+,-\}\),
\begin{align*}
  \partial_{k}^{\epsilon}((x:\obj)\vdashps x:\obj)
  &= (x:\obj)
  &   \nameref{rule:pss}
  \\
  \partial_{k}^{\epsilon}(\Gamma,y:A,f:\arr[A]{x}{y}\vdashps f\arr[A]{x}{y})
  & =
    \begin{cases}
      \partial_{k}^{\epsilon}(\Gamma\vdashps x:A)
      & \dim A > k\\
      (\partial_{k}^{\epsilon}(\Gamma\vdashps x:A),y:A,f:\arr[A]{x}{y})
      & \dim A \leq k-1\\
      \partial_{k}^{-}(\Gamma\vdashps x:A)
      & \dim A = k,\epsilon = -\\
      (\rmlast\partial_{k}^{+}(\Gamma\vdashps x:A),y:A)
      & \dim A = k, \epsilon = +
    \end{cases}
  &
    \nameref{rule:pse}
  \\
  \partial_{k}^{\epsilon}(\Gamma\vdashps y:A)
  &=  \partial_{k}^{\epsilon}(\Gamma\vdashps f : \arr[A]{x}{y})
  &  \nameref{rule:psd}\\
  \partial_{k}^{\epsilon}(\Gamma\vdashps)
  &= \partial_{k}^{\epsilon}
    (\Gamma\vdashps x:\obj) & \nameref{rule:ps}
\end{align*}
The first author proved~\cite{benjamin_formalisation_2022} that
\(\partial_{k}^{\pm}\Gamma\) are actually pasting contexts. By weakening, we
get two substitutions \(\Gamma\vdash s_{k}^{\Gamma} : \partial_{k}^{-}\Gamma\)
and \(\Gamma \vdash t_{k}^{\Gamma} : \partial_{k}^{+}\Gamma\) that assign each
variable to the variable with the same name.
Fig.~\ref{fig:rules-catt} then completely defines the dependent type theory
\catt. The definition given here differs slightly from the original definition
due to Finster and Mimram~\cite{finster_typetheoretical_2017}. We unite the two
rules they present into a single one and we use an alternative, but equivalent,
second side condition. The equivalence has been shown by the first
author~\cite[Section~3.5.1]{benjamin_type_2020}. We note that the two side
conditions can be further unified into a single one
\begin{align*}
  \Var(s^{\Gamma}_{\dim A}) &= \Var(u)\cup \Var(A) &
  \Var(t^{\Gamma}_{\dim A}) &= \Var(v)\cup \Var(A)
\end{align*}
using that \(s^{\Gamma}_k  = t^{\Gamma}_k = \id\) for \(k\ge \dim \Gamma\) and
a result of the first author~\cite[Lemma~76]{benjamin_type_2020}. Moreover,
we observe that the side condition can only be used when
\(\dim \Gamma \le \dim A + 1\).

\section{Computads for weak \texorpdfstring{\(\omega\)}{ω}-categories}
\label{sec:computads}

Computads are combinatorial structures for presenting higher categories, first
invented by Street~\cite{street_limits_1976}, and later generalised by
Burroni~\cite{burroni_higherdimensional_1993} and
Batanin~\cite{batanin_computads_1998}. They consist of sets of generators with
specified source and target cells, which are cells of the higher category freely
generated by the computad. The apparent circularity of this description of
computads is resolved by defining the category of computads and free
higher categories on computads mutually recursively.

Following Dean, Finster, Markakis, Reutter and
Vicary~\cite{dean_computads_2022}, to define computads for
\(\omega\)\-categories, we need to start by defining globular sets and
Batanin trees. A \emph{globular set} \(X\) consists of a set
\(X_n\) for every \(n\in\N\) together with source and target functions
\( \src, \tgt \colon X_{n+1} \to X_n\) subject to the globularity conditions
\begin{align*}
  \src \circ \src &= \src \circ \tgt &
  \tgt \circ \src &= \tgt \circ \tgt
\end{align*}
We call elements of \(X_n\) the \emph{\(n\)\-cells} of \(X\), and we call
\(\src\) and \(\tgt\) the \emph{source} and \emph{target} functions. Morphisms
of globular sets \(f \colon X\to Y\) are sequences of functions
\(f_n\colon X_n\to Y_n\) commuting with the source and target functions.
Globular sets with their morphisms form a category \(\gSet\), which is a
presheaf topos. We denote by \(\disk^n\) the globular set represented by
\(n\in\N\), and we denote by \(\sphere^{n-1}\) the globular set obtained by
removing the unique \(n\)\-cell of \(\disk^n\). We will also denote the
obvious inclusion by \(\iota_n \colon \sphere^{n-1} \to \disk^n\).

Batanin trees are a family of globular sets, corresponding to \emph{pasting
schemes}, that play the role of arities for the operations of
\(\omega\)\-categories. The set of Batanin trees \(\Bat\) is defined inductively
by one constructor
\[
  \br\colon \List \Bat\to \Bat.
\]
In particular, there exists a Batanin tree \(\br[]\) with no branches. The
\emph{dimension} of a Batanin tree is the height of the corresponding tree, and
it is defined recursively by
\begin{align*}
  \dim (\br[])
    &= 0 &
  \dim (\br[B_1,\dots,B_n])
    &= \max(\dim B_1 + 1, \dim (\br[B_2,\dots,B_n])).
\end{align*}
The
\emph{\(k\)\-boundary} of a Batanin tree \(B\) is the Batanin tree defined
recursively by
\begin{align*}
  \partial_0 B
    &= \br[] &
  \partial_{k+1} \br[B_1,\dots,B_n]
    &= \br[\partial_k B_1,\dots,\partial_k B_n].
\end{align*}
A globular set \(\Pos{B}\) can be assigned to each Batanin tree
\(B\), whose cells we call \emph{positions}. This globular set can be generated
inductively, as shown by Dean et al.~\cite{dean_computads_2022}. Equivalently,
it can be described recursively by the formula
\[
  \Pos(\br[B_1,\dots,B_n]) = \bigvee_{i=1}^n \Sigma\Pos(B_i).
\]
as shown in our previous work~\cite{benjamin_opposites_2024}.
Here the suspension operation \(\Sigma\) sends a globular set \(X\) to the
bipointed globular set with two \(0\)\-cells \(v_-\) and \(v_+\) and globular
set of cells from \(v_-\) to \(v_+\) given by \(X\), and the wedge sum of two
bipointed globular sets \(X \vee Y\) is obtained by the coproduct of \(X\) and
\(Y\) by identifying the second basepoint of \(X\) with the first basepoint of
\(Y\). Positions of a Batanin tree are
precisely the sectors of the corresponding planar tree, as explained by
Berger~\cite{berger_cellular_2002}. We illustrate this correspondence in
Fig.~\ref{fig:ctx-trees-ps}. The positions of the \(k\)\-boundary of \(B\) may
be included back in the positions of \(B\) via the source and target inclusions
\[
  s^B_k, t_k^B \colon \Pos(\partial_k B) \to \Pos (B)
\]
which are similarly defined recursively by \(s_0^B\) and \(t_0^B\) being
the morphisms out of the \(0\)\-disk picking the left and right basepoint of
\(\Pos(B)\) and
\begin{align*}
  s_{k+1}^{\br[B_1,\dots,B_n]} &= \bigvee_{i=1}^n \Sigma s_k^{B_i} &
  t_{k+1}^{\br[B_1,\dots,B_n]} &= \bigvee_{i=1}^n \Sigma t_k^{B_i}.
\end{align*}
We note that the source and target inclusions, satisfy equations dual to the
globularity conditions, and that they are identities for \(k \ge \dim B\).

The category of \(n\)\-computads \(\Comp_n\) is then defined inductively on
\(n\in\N\) together with four functors
\begin{align*}
  \Free_n &\colon \gSet \to \Comp_n &
  \Cell_n &\colon \Comp_n \to \Set \\
  \tr_{n-1} &\colon \Comp_{n} \to \Comp_{n-1} &
  \Sphere_n &\colon \Comp_n \to \Set.
\end{align*}
The functor \(\Free_{n}\) views a globular set as an \(n\)\-computad. The
functor \(\tr_{n-1}\) forgets the top-dimensional generators of a computad. The
functors \(\Cell_{n}\) and \(\Sphere_{n}\) send an \(n\)\-computad to the set of
\(n\)\-cells of the \(\omega\)\-category it presents, and the set of pairs of
parallel \(n\)\-cells respectively. In the same mutual induction, three natural
transformations are defined
\begin{align*}
  \pr_1 &\colon \Sphere_n \Rightarrow \Cell_n &
  \pr_2 &\colon \Sphere_n \Rightarrow \Cell_n &
  \bdry &\colon \Cell_n \Rightarrow \Sphere_{n-1}\tr_{n-1}
\end{align*}
where \(\pr_i\) are the obvious projections and \(\bdry\) takes a cell to its
source and target. Finally, for every Batanin tree \(B\), a subset of spheres
\(\Full_n(B)\subseteq \Sphere_n(\Free_n\Pos(B))\) is defined to be \emph{full}.

The induction starts by letting \(\Comp_{-1}\) be the terminal category and
\(\Sphere_{-1}\) the functor picking the terminal set. For \(n\in\N\), an
\(n\)\-computad \(C\) consists of an \({(n-1)}\)\-computad \(C_{n-1}\), a set
of generators \(V_n^C\), and an attaching function
\(\phi_n^C\colon V_n^C \to \Sphere_{n-1}(C_{n-1})\). A
morphism \(f \colon C\to D\) of \(n\)\-computads consists of a morphism
\(f_{n-1} \colon C_{n-1}\to D_{n-1}\) and a function
\(f_{V} \colon V_n^C \to \Cell_n(D)\) such that
\[
  \bdry_{n,D} \circ f_{n,V} = \Sphere_{n-1}(f_{n-1}) \circ \phi_n^C.
\]
The functor \(\tr_{n-1}\) is the obvious first projection.

The set of \(n\)\-cells of a computad \(C\) is defined inductively by two
constructors. The first constructor \(\var\) produces an \(n\)\-cell out of
a generator \(v\in V_n^C\). The second constructor \(\coh\) is only used when
\(n > 0\) and produces an \(n\)\-cell out of a Batanin tree \(B\) of dimension
at most \(n\), a full \((n-1)\)\-sphere \(A\) of \(B\) and a morphism
\(f \colon \Free_n\Pos(B)\to C\). The boundaries of those cells are given
recursively by
\begin{align*}
  \bdry_{n,C}(\var v) &= \phi_n^C(v) &
  \bdry_{n,C}(\coh(B,A,f)) &= \Sphere_{n-1}(f_{n-1})(A)
\end{align*}
An \(n\)\-sphere of an \(n\)\-computad \(C\) is then a pair of \(n\)\-cells
\((a,b)\) with the same boundary, and the projection natural transformations are
the obvious ones. Composition of morphisms and the action of
a morphism on \(n\)\-cells are given mutually recursively by
\begin{gather*}
  \Cell_n(g)(\var v) = g_{V}(v) \\
  \Cell_n(g)(\coh(B,A,f)) = \coh(B,A,g\circ f) \\
  g\circ f = (g_{n-1}\circ f_{n-1}, \Cell_n(g)\circ f_{n,V})
\end{gather*}
The \(n\)\-computad \(\Free_n X\) consists of \(\Free_{n-1} X\), the set
\(X_n\) of \(n\)\-cells and the attaching function
\[
  \phi^{\Free X}_n (x) = (\var \src x, \var \tgt x).
\]
The fullness condition for a sphere \((a,b)\in\Sphere_n(\Free_n \Pos(B))\) is
a condition on the \emph{supports} of the cells \(a\) and \(b\), which are the
sets of \(n\)\-positions used in the definition of \(a\) and \(b\) respectively.
More precisely, the support of a cell \(c\in \Cell_n(C)\) is the set of
\(n\)\-dimensional generators of \(C\) defined by
\begin{align*}
  \supp_n(\var v) &= \{v\} &
  \supp_n(\coh(B,A,f)) &= \bigcup_{p\in \Pos_n(B)} \supp_n(f_{V}(p)).
\end{align*}
A sphere \((a,b)\) is then declared to be full when the supports of \(a\) and
\(b\) contain precisely the \(n\)\-positions in the image of the source and
target inclusions \(s_n^B\) and \(t_n^B\) respectively, and the sphere
\(\bdry_n(a)\) is full.

Finally, the category of computads \(\Comp\) is the limit of the categories of
\(n\)\-computads over the obvious truncation functors \(\tr_n\). We will denote
the projections out of the limit by \(\Tr_n \colon \Comp \to \Comp_n\). By the
universal property of the limit, the functors \(\Free_n\) give rise to a
functor \[\Free \colon \gSet\to \Comp,\]
while the functor \(\Cell_n\Tr_n\) together with the natural transformations
\(\pr_i\bdry_n\) give rise to a functor in the opposite direction
\[\Cell \colon \Comp\to \gSet.\]
It is shown by Dean et al.~\cite{dean_computads_2022} that the functor \(\Free\)
is left adjoint to \(\Cell\); the unit being the morphism of globular sets given
by \(\var\) and the counit being the morphism of computads viewing cells of
\(C\) as generators of \(\Free\Cell C\). The algebras of the monad
\(T\) induced by this adjunction are precisely the \(\omega\)\-categories of
Leinster~\cite{leinster_operads_2000}.

\subsection{The category with families structures on
  \texorpdfstring{\(\gSet\)}{gSet} and \texorpdfstring{\(\Comp\)}{Comp}}

The opposite of the category of globular sets \(\gSet^{\op}\) can be equipped
with the structure of a category with families. The presheaf of families
\(T^{\gSet} \colon \gSet\to \Fam\) sends an object \(X\) to the family
\begin{align*}
  \Type^{\gSet}(X) &= \coprod\nolimits_{n\in\N}\gSet(\sphere^{n-1},X) &
  \Term^{\gSet}(X;(n,A)) &= \left\{ x \colon \disk^n \to X \;\middle|\;
    x\circ \iota_n = A\right\}
\end{align*}
By the Yoneda lemma and the definition of \(\sphere^n\), we see that types of
a globular set \(X\) are in bijection to pairs of cells of \(X\) with common
source and target; terms of a type \((n,A)\) are \(n\)\-cells of \(X\) whose
source and target are given by \(A\). The extension of a globular set \(X\) by a
type \((n,A)\) is given by the following pushout square in \(\gSet\),
\[\begin{tikzcd}[column sep = large]
  \sphere^{n-1}\ar[r, "\iota_{n}", hook]\ar[d,"A"']
  \ar[dr,phantom,"\ulcorner"very near end]
  & \disk^{n} \ar[d,"v_{X,A}"] \\
  X\ar[r,"\pi_{X,A}"'] & X.A
\end{tikzcd}\]
or more concretely by adjoining a new cell of dimension \(n\) to \(X\) whose
source and target are determined by \(A\). The universal property of the
extension \(X.A\) in the category with families structure is exactly the
universal property of the pushout in the category \(\gSet\).

The opposite of the category of computads \(\Comp^{\op}\) can similarly equipped
with the structure of a category with families by sending a computad \(C\) to
the family
\begin{align*}
  \Type^{\Comp}(C) &= \coprod\limits_{n\in\N} \Comp(\Free\sphere^{n-1},C) &
  \Term^{\Comp}(C;(n,A)) &= \left\{ c\colon \Free\disk^n \to C \;\middle|\;
    c\circ \Free\iota_n = A\right\}
\end{align*}
By representability of \(\Sphere_n\) and
\(\Cell_n\)~\cite[Corollary~4.2]{dean_computads_2022}, we have that the set of
types of a computad \(C\) is the union of the sets of \(n\)\-spheres of \(C\)
over all \(n\), and that the set of terms of a type \((n,A)\) is the set of
cells with boundary sphere \(A\). The extension of a computad \(C\) by a type
\((n,A)\) is given similarly by the following pushout square in \(\Comp\).
\[
  \begin{tikzcd}[column sep = large]
    \Free\sphere^{n}\ar[r, "\Free\iota_{n+1}", hook]\ar[d,"A"']
    \ar[dr,phantom,"\ulcorner"very near end]
    & \Free\disk^{n+1} \ar[d,"v_{C,A}"] \\
    C\ar[r,"\pi_{C,A}"'] & C.A
  \end{tikzcd}
\]
which exists by~\cite[Proposition~6.3]{dean_computads_2022}; it is obtained by
adjoining a new generator of dimension \(n\) to \(C\) with boundary sphere
\(A\).
By construction, the functor \(\Free \colon \gSet\to \Comp\) preserves the
chosen pushout square, so the functor
\(\Free^{\op}\colon\gSet^{\op}\to \Comp^{\op}\) is part of a morphism of
categories with families, equipped with the natural transformations
given by the action of \(\Free\) on morphisms.

\section{Globular sets and GSeTT}
\label{sec:globular-gsett}

The type theory \gsett has been studied in detail by Benjamin, Finster and
Mimram~\cite{benjamin_globular_2021} where it was shown that its category of
\(\Set\)-models is equivalent to the category of globular sets, and that its
syntactic category is equivalent to the opposite of the subcategory of finite
globular sets. Here by finite globular sets, we mean globular sets \(X\) such
that \(\amalg_{n\in\N} X_n\) is finite. More explicitly, using that contexts,
types, terms and substitutions of \gsett have unique derivations, we can define
a morphism of categories with families \(V\colon \syn_{\gsett}\to \gSet^{\op}\)
recursively on the syntax by the rules given in
Fig.~\ref{fig:comparison-functor-gsett}. Functoriality and naturality of the
assignment \(V\) are straightforward to check by mutual induction, while
preservation of the extension operation is immediate from the definition.

\begin{figure}
  \centering
  \fbox{
    \begin{minipage}{.95\linewidth}
      \begin{mathpar}
        V(\emptyset) = \emptyset \and
        V_{\Gamma, \emptyset}(\sub{}) = \sub{} \and V_\Gamma^{\Type}(\obj) = \sub{} \\
        \begin{tikzcd}[ampersand replacement=\&]
          \sphere^{\dim A} \& \disk^{\dim A + 1} \\
          V(\Gamma) \& V(\Gamma, x:A)
          \ar[from=1-1, to=1-2, hook, "\iota_{\dim A}"]
          \ar[from=2-1, to=2-2, dashed, "\pi_{\Gamma, A}"']
          \ar[from=1-1, to=2-1, "V^{\Type}_\Gamma(A)"']
          \ar[from=1-2, to=2-2, dashed, "v_{\Gamma,A}"]
          \ar[from=1-1, to=2-2, phantom, "\ulcorner"{very near end}]
        \end{tikzcd} \and
        \begin{aligned}
          V_{\Delta,(\Gamma,x: A)} (\sub{\gamma,t} )
          &= \sub{V_{\Delta,\Gamma}(\gamma),V^{\Term}_{\Delta,A[\gamma]}(t)} \\
          V_{\Gamma}^{\Type}(\arr[A]{u}{v})
          &= \sub{V_{\Gamma,A}^{\Term}(u), V_{\Gamma,A}^{\Term}(v)} \\
          V^{\Term}_{(\Gamma, x:A), B}(\var v)
          &=
            \begin{cases}
              v_{\Gamma,A} & \text{if } v=x \\
              \pi_{\Gamma, A} \circ V^{\Term}_{\Gamma,B}(\var v) &\text{otherwise}
            \end{cases}
        \end{aligned}
      \end{mathpar}
    \end{minipage}
  }
  \caption{The morphism \(V\colon \syn_{\gsett}\to \gSet^{\op}\)}
  \label{fig:comparison-functor-gsett}
\end{figure}

Since the type theory \gsett has no term constructors, the cells of the
globular set \(V\Gamma\) are precisely the valid terms of the context
\(\Gamma\), that is, its variables. It follows that the functor constructed here is isomorphic
to the functor introduced by Benjamin, Finster and
Mimram~\cite[Definition~15]{benjamin_globular_2021}, hence fully faithful with
essential image the finite globular sets. Fullness and faithfulness amount to
morphisms of globular sets \(V\Gamma\to V\Delta\) being functions from the
variables of \(\Delta\) to the variables of \(\Gamma\) compatible with the
source and target functions. Essential surjectivity onto finite globular sets
amounts to every finite globular set \(X\) being a finite extension
\(\emptyset.A_1.A_2.\cdots.A_n\) for some finite sequence of spheres
\(A_k\in \Type^{\gSet}(\emptyset.A_1.A_2.\cdots.A_{k-1})\). For similar reasons,
the natural transformations \(V^{\Type}\) and \(V^{\Term}\) are invertible.

\section{Pasting contexts and Batanin trees}
\label{sec:pasting-batanin}

Before comparing the type theory \catt with the category of computads, we need
to compare the arities of the operation \(\coh\) in the two theories, namely the
set of pasting contexts \(\PsCtx\) and the set of Batanin trees \(\Bat\). To
compare the two notions, we will use the auxiliary notions of \emph{smooth
zigzag sequence} and Street's \emph{globular cardinal}~\cite{street_petit_2000},
and their boundaries. The correspondence between the four notions is illustrated
in Fig.~\ref{fig:ctx-trees-ps} for a particular globular cardinal.

Following Weber~\cite[Section~4]{weber_generic_2004}, we may assign functorially
to every globular set \(X\), a preorder \(\Sol(X)\), whose carrier is the set
\(\amalg_{n\in\N}X_n\) of cells of \(X\) and whose relation
\(\blacktriangleleft\) is the reflexive and transitive closure of the relation
\(\prec\) defined by \(\src x \prec x\) and \(x \prec \tgt x\) for every cell
\(x\). A \emph{globular cardinal} is a globular set \(X\) for which
\(\Sol X\) is a non-empty finite total order. We define further the
\(k\)\-boundary of a globular cardinal \(X\) to be the globular cardinal
\(\partial_k X\) whose \(n\)\-cells are given by
\[
  (\partial_k X)_n = \begin{cases}
    X_n & \text{if } n < k \\
    (X_k)_{/\sim} & \text{if } n = k \\
    \emptyset & \text{if } n > k
  \end{cases}
\]
where two \(k\)\-cells are related by \(\sim\) if and only if they have the same
source and target. The source and target inclusions
\(s_k^X, t_k^X \colon \partial_k X\to X\) are the morphisms of globular sets
that are identity on cells of dimension less than \(k\), and they pick the
least and greatest representative of each equivalence class of \(k\)\-cells
respectively. It is easy to see that the \(k\)\-boundary operation is
functorial and the source and target inclusions are natural with respect to
isomorphisms. The notion of
boundary we define here is a mild generalisation of
the boundary of a globular cardinal as defined by Benjamin, Finster and
Mimram~\cite{benjamin_globular_2021}, where they define the source and target
\(\partial^\pm X\) of a globular cardinal instead. The precise relation between
the two notions is that
\begin{align*}
  \partial^+ X = \image(s_{\dim X -1}^X) &&
  \partial^- X = \image(t_{\dim X -1}^X).
\end{align*}
where \(\dim X\) is the maximum of the dimension of the cells of \(X\).

\begin{lemma}\label{lem:susp-card}
  The suspension \(\Sigma X\) of a globular cardinal \(X\) is a globular
  cardinal whose basepoints are its minimum and maximum cell. Moreover, the
  following equalities hold:
  \begin{align*}
    \partial_{k+1} (\Sigma X) &= \Sigma (\partial_k X) &
    s_{k+1}^{\Sigma X} &= \Sigma s_k^X &
    t_{k+1}^{\Sigma X} &= \Sigma t_k^X.
  \end{align*}
\end{lemma}
\begin{proof}
  The total order \(\Sol(\Sigma X)\) is obtained from the order \(\Sol X\) by
  freely adjoining a minimum and maximum element. The equalities follow easily
  from the definition of the suspension.
\end{proof}
\begin{lemma}\label{lem:wedge-card}
  The wedge sum \(X\vee Y\) of two cardinals \(X\) and \(Y\) whose basepoints
  are their minimum and maximum cells is again a globular cardinal with the same
  property. Moreover, there exists an isomorphism, making the following two
  triangles commute:
  \begin{align*}
    \begin{tikzcd}[ampersand replacement=\&, column sep = 0pt, row sep = small]
      \partial_{k+1}(X\vee Y) \&\&
      (\partial_{k+1}X) \vee (\partial_{k+1}Y) \\
      \& X\vee Y
      \ar[from=1-1, to=1-3, "\sim"]
      \ar[from=1-1, to=2-2, "s_{k+1}^{X\vee Y}"']
      \ar[from=1-3, to=2-2, "s_{k+1}^{X}\vee s_{k+1}^{Y}"]
    \end{tikzcd} &&
    \begin{tikzcd}[ampersand replacement=\&, column sep = 0pt, row sep = small]
      \partial_{k+1}(X\vee Y) \&\&
      (\partial_{k+1}X) \vee (\partial_{k+1}Y) \\
      \& X\vee Y
      \ar[from=1-1, to=1-3, "\sim"]
      \ar[from=1-1, to=2-2, "t_{k+1}^{X\vee Y}"']
      \ar[from=1-3, to=2-2, "t_{k+1}^{X}\vee t_{k+1}^{Y}"]
    \end{tikzcd}
  \end{align*}
\end{lemma}
\begin{proof}
  The order \(\Sol(X\vee Y)\) is obtained from by identifying the
  maximum cell of \(X\) with the minimum cell of \(Y\). The isomorphism is the
  identity on cells of dimension at most \(k\), and it is given on
  \((k+1)\)\-cells by the isomorphism of sets
  \[
    (X_{k+1} \amalg Y_{k+1})_{/\sim} \xrightarrow{\sim}
    (X_{k+1})_{/\sim} \amalg (Y_{k+1})_{/\sim}
  \]
  which exists, since no element of \(X_{k+1}\) is related by \(\sim\) to
  an element of \(Y_{k+1}\). Commutativity of the triangles follows easily by
  the definition of the source and target inclusions.
\end{proof}

\begin{proposition}\label{prop:batanin-card}
  The globular set \(\Pos(B)\) is a globular cardinal for every Batanin tree
  \(B\). Moreover, there exists an isomorphism
  \(\partial_k\Pos(B)\cong\Pos(\partial_k B)\) making the following diagrams
  commute
  \begin{align*}
    \begin{tikzcd}[ampersand replacement=\&, column sep = 0pt, row sep = small]
      \partial_k\Pos(B) \&\& \Pos(\partial_k B) \\
      \&\Pos(B)
      \ar[from=1-1, to=1-3, "\sim"]
      \ar[from=1-1, to=2-2, "s_k^{\Pos(B)}"']
      \ar[from=1-3, to=2-2, "s_k^{B}"]
    \end{tikzcd} &&
    \begin{tikzcd}[ampersand replacement=\&, column sep = 0pt, row sep = small]
      \partial_k\Pos(B) \&\& \Pos(\partial_k B) \\
      \&\Pos(B)
      \ar[from=1-1, to=1-3, "\sim"]
      \ar[from=1-1, to=2-2, "t_k^{\Pos(B)}"']
      \ar[from=1-3, to=2-2, "t_k^{B}"]
    \end{tikzcd}
  \end{align*}
\end{proposition}
\begin{proof}
  This is an immediate consequence of Lemmas~\ref{lem:susp-card}
  and~\ref{lem:wedge-card}, and the fact that the unit of the wedge sum
  \(\disk^0\) is also a globular cardinal.
\end{proof}

We now show every globular cardinal is isomorphic to one of the form
\(\Pos(B)\) for unique \(B\). To do so, we use the
\emph{smooth zigzag sequences} introduced by Weber~\cite{weber_generic_2004}.
A smooth exact sequence is a non-empty finite sequence of natural numbers
\((m_1,\dots,m_n)\) such that
\(m_1 = m_n = 0\) and \(\left|m_i - m_{i-1}\right| = 1\) for all \(i < n\).
To each globular cardinal \(X\), we may assign a smooth zigzag sequence
\(\Zig X = (\dim x_1,\dots,\dim x_k)\), where \(x_1\prec \dots \prec x_k\) are
the cells of \(X\) in increasing order. Conversely, from a smooth zigzag
sequence \((m_1,\dots,m_n)\), we may construct a globular cardinal
\(\Card(m_1,\dots,m_n)\), whose \(k\)\-cells are indices \(i\le n\) such that
\(m_i = k\) and whose source and target functions are given by
\begin{align*}
  \src(i) &= \max\{j < i \ |\ m_j = m_i - 1\} &
  \tgt(i) &= \min\{j > j \ |\ m_j = m_i - 1\}
\end{align*}
Isomorphic globular cardinals are sent to the same smooth zigzag sequence, and
\begin{align*}
  \Zig \circ \Card (m_\bullet) &= m_\bullet &
  \Card \circ \Zig (X) &\cong X &
\end{align*}
so that smooth zigzag sequences are in bijection with isomorphism classes of
globular cardinals.

To show that \(\Pos\) is a bijection onto isomorphism classes of globular
cardinals, it suffices to show that \(\Zig\circ \Pos\) is a bijection onto
smooth zigzag sequences. From the proofs of Proposition~\ref{prop:batanin-card},
we can derive a recursive formula for the composite \(\Zig\circ\Pos\):
\[
  \Zig \Pos(\br[B_1,\dots,B_n]) = \bigHash_{i=1}^n (\Zig \Pos(B_1))^+
\]
where \((m_1,\dots,m_n)^+ = (0,m_1+1,\dots,m_n+1,0)\) and \(\#\) is the
operation concatenating two smooth zigzag sequences by identifying the last
element of the first with the first element of the second. The inverse of
\(\Zig\circ \Pos\) is given by the function \(\Tree\) defined recursively on the
length of a sequence by
\begin{gather*}
  \Tree(0,m^1_1,\dots,m^1_{n_1},0,\dots,0,m^l_1,\dots,m^l_{n_l},0)
    = \br[B_1,\dots,B_l] \\
  B_i = \Tree(m^i_1-1,\dots,m^i_{n_i}-1)
\end{gather*}
where \(m_i^j > 0\) for all indices \(i, j\). Hence, we have proven the
following:

\begin{proposition}\label{prop:card-batanin}
  For each globular cardinal \(X\), there exists unique Batanin tree \(B\)
  such that \(\Pos(B)\cong X\).
\end{proposition}

The last step in comparing pasting contexts to Batanin trees is to show that
the former are also in bijection to isomorphism classes of globular cardinals.
This has already been shown by Benjamin, Finster and Mimram. The globular set
\(V\Gamma\) associated to a
pasting context \(\Gamma\), seen as a context of \gsett, is a globular
 cardinal~\cite[Proposition~40]{benjamin_globular_2021} and conversely every
globular cardinal is isomorphic to \(V\Gamma\) for a unique pasting context
\(\Gamma\)~\cite[Proposition~42]{benjamin_globular_2021}. Intuitively, the
order \(\Sol(V\Gamma)\) is the order in which the variables appear in the
right hand side of the auxiliary judgements in the derivation tree of
\(\Gamma\vdashps\). Conversely, to each smooth zigzag sequence
\((m_0,\dots,m_k)\), we may assign a derivation of the form \(\Gamma\vdashps\)
by starting from the rule \nameref{rule:pss}, iteratively using the rule
\nameref{rule:pse} when \(m_i > m_{i-1}\) and the rule \nameref{rule:psd}
otherwise, until the sequence is exhausted, and concluding with the rule
\nameref{rule:ps}. Since the derivation of the judgement \(\Gamma\vdashps\) is
unique, one can check that \(\Zig\circ V\) is a bijection with inverse the
assignment described above. The following Theorem summarises
the results of this section.

\begin{theorem}\label{thm:psctx-trees}
  There exists a bijection \(B_\bullet \colon \PsCtx\to \Bat\) together with a
  family of isomorphisms
  \[R^{\Bat}_\Gamma\colon V\Gamma\xrightarrow{\sim}\Pos(B_\Gamma)\]
  sending the image of
  \(V(s^{\Gamma}_{\dim \Gamma-1})\) to the image of
  \(s_{\dim \Gamma-1}^{B_\Gamma}\), and similarly for the target inclusions.
\end{theorem}
\begin{proof}
  The bijection \(B_\bullet\) is the composite \(\Tree\circ\Zig\circ V\). By
  Proposition~\ref{prop:batanin-card} and naturality, the following diagram
  commutes:
  \[\begin{tikzcd}[ampersand replacement=\&]
    {\partial_{k} (V\Gamma)} \& {\partial_k (\Pos(B_\Gamma))} \& {\Pos(\partial_k B_\Gamma)} \\
    V\Gamma \& {\Pos(B_\Gamma)}
    \arrow["{\partial R^{\Bat}_\Gamma}", from=1-1, to=1-2]
    \arrow["{s^{V\Gamma}_k}"', from=1-1, to=2-1]
    \arrow["{R^{\Bat}_\Gamma}"', from=2-1, to=2-2]
    \arrow["{s^{\Pos(B_\Gamma)}_k}"', from=1-2, to=2-2]
    \arrow["\sim", from=1-2, to=1-3]
    \arrow["{s^{B_\Gamma}_k}", from=1-3, to=2-2]
  \end{tikzcd}\]
  Hence, the isomorphism \(R_\Gamma^{\Bat}\) identifies the image of
  \(s_k^{V\Gamma}\) with that of \(s_k^{B_\Gamma}\). The former for
  \(k = \dim \Gamma -1\) has been shown to be equal to the image of
  \(V(s^\Gamma_k)\) by Benjamin, Finster and
  Mimram~\cite[Lemma~45]{benjamin_globular_2021}.
\end{proof}

\begin{figure}
  \centering
  \begin{tabular}{c@{\qquad}c@{\qquad}c@{\qquad}c}
    \(\scriptstyle  x\prec f \prec \alpha \prec g \prec y \prec h \prec z\) &
    \(\scriptstyle  (0,1,2,1,0,1,0)\) &
    \(\scriptstyle  \br[\br[\br[]],\br[]]\) &
    \(\scriptstyle \nameref{rule:pss}\nameref{rule:pse}^2\nameref{rule:psd}^2
      \nameref{rule:pse}\nameref{rule:psd}\nameref{rule:ps}\)\\
    \(\begin{aligned}\begin{tikzcd}[ampersand replacement=\&]
      x
      \ar[r, "f"{name=0}, bend left = 45]
      \ar[r, "g"{name=1, swap}, bend right = 45] \&
      y \ar[r, "h"] \&
      z
      \ar[from = 0, to = 1, Rightarrow, "\alpha",
        shorten <= 5pt, shorten >= 5pt]
    \end{tikzcd}\end{aligned}\) &
    \(\begin{aligned}\begin{tikzpicture}[xscale=0.6,yscale=0.7, inner sep=0]
      \draw[-,dotted] (0,1) -- (6.5,1);
      \draw[-,dotted] (0,2) -- (6.5,2);
      \draw[->,gray] (0,0) -- (6.5,0);
      \draw[->,gray] (0,0) -- (0,2.5);
      \node[label={[label distance = 3]85:{$\scriptstyle x$}}] (0) at (0,0) {$\bullet$};
      \node[label={[label distance = 3]above:{$\scriptstyle f$}}] (1) at (1,1) {$\bullet$};
      \node[label={[label distance = 3]above:{$\scriptstyle \alpha$}}] (2) at (2,2) {$\bullet$};
      \node[label={[label distance = 3]above:{$\scriptstyle g$}}] (3) at (3,1) {$\bullet$};
      \node[label={[label distance = 3]above:{$\scriptstyle y$}}] (4) at (4,0) {$\bullet$};
      \node[label={[label distance = 3]above:{$\scriptstyle h$}}] (5) at (5,1) {$\bullet$};
      \node[label={[label distance = 3]above:{$\scriptstyle z$}}] (6) at (6,0) {$\bullet$};
      \draw[-] (0) -- (1) -- (2) -- (3) -- (4) -- (5) -- (6);
    \end{tikzpicture}\end{aligned}\) &
    \(\begin{aligned}
      \begin{tikzcd}[ampersand replacement=\&, row sep = small,
        column sep = tiny]
        \overset{\alpha}{\bullet} \\
        \bullet \& {}\& \overset{h}{\bullet} \\
        \& \bullet
        \arrow[no head, from=1-1, to=2-1, shorten <= -3pt, shorten >= -3pt,
          "g"{very near end}, "f"{very near end, swap}]
        \arrow[no head, from=2-1, to=3-2, shorten <= -3pt, shorten >= -3pt,
          "x"{very near end, swap}]
        \arrow[no head, from=2-3, to=3-2, shorten <= -3pt, shorten >= -3pt,
          "z"{very near end}]
        \arrow[from=3-2, to=2-2, phantom, "\scriptstyle y"{very near start}]
      \end{tikzcd}\end{aligned}\) &
      \(\begin{aligned}
        (x : \obj)(y : \obj)(f : \arr{x}{y}) \\
        (g : \arr{x}{y})(\alpha : \arr{f}{g}) \\
        (z : \obj)(h : \arr{y}{z}) \\
      \end{aligned}\) \\
      \scriptsize{(a) A globular cardinal} &
      \scriptsize{(b) A smooth zigzag sequence} &
      \scriptsize{(c) A Batanin tree} &
      \scriptsize{(d) A pasting context}
  \end{tabular}
  \caption{The correspondence between globular cardinals, smooth zigzag sequences, Batanin trees and pasting contexts.}
  \label{fig:ctx-trees-ps}
\end{figure}

\section{Computads and CaTT}
\label{sec:main}

We are now in the position to compare the type theory \catt with the category of
computads, extending the functor \(V\) of Section~\ref{sec:globular-gsett} by
using the bijection \(B_\bullet\) of Section~\ref{sec:pasting-batanin}. We
define a morphism of categories with families \(R \colon \syn_{\catt}\to
\Comp^{\op}\) by induction on the syntax, as indicated in
Fig.~\ref{fig:comparison-functor-catt}. Mutually inductively, we verify that the
following diagram commutes
\[\begin{tikzcd}
  \syn_{\gsett} & \gSet^{\op} \\
  \syn_{\catt} & \Comp^{\op}
  \arrow["V", from=1-1, to=1-2]
  \arrow["R"', from=2-1, to=2-2]
  \arrow["J"', hook, from=1-1, to=2-1]
  \arrow["\Free", from=1-2, to=2-2]
\end{tikzcd}\]
where the functor \(J\) views contexts and substitutions of \gsett as contexts
and substitutions of \catt. For the last case in
Fig.~\ref{fig:comparison-functor-catt}, we implicitly use the fact that
morphisms \(\sphere^n\to C\) are in natural bijection with elements of
\(\Sphere_n(C)\). Moreover, for this case to be well-defined, we need to show
that for a type \(A\) satisfying the side condition of the
rule~\nameref{rule:coh}, the sphere
\((\Free R_\Gamma^{\Bat}\circ R^{\Type}A)\) is again full. Commutativity of the
square is immediate by \(\Free\) preserving strictly the chosen pushouts of
sphere inclusions. To prove fullness, we see from the inductive
presentation that we may canonically identify variables of \(\Gamma\) with the
generators of \(R\Gamma\). Under this identification, it suffices to prove that
\[
  \supp(R^{\Term}_{\Gamma,A} u) = \Var(u) \cup \Var(A)
\]
Indeed, for a type \(A' =\arr[A]{u}{v}\) satisfying the side condition, we have
by the support lemma~\cite[Lemma~7.3]{dean_computads_2022} and
Theorem~\ref{thm:psctx-trees}, we see that
\[\begin{split}
  \supp(\Free R_\Gamma^{\Bat} \circ R^{\Term}_{\Gamma,A}(u))
    &= R_\Gamma^{\Bat}(\supp R^{\Term}_{\Gamma,A}(u))
    = R_\Gamma^{\Bat}(\Var(u) \cup \Var(A)) \\
    & = R_\Gamma^{\Bat}(\Var(s_{\dim A'}^\Gamma))
    = \image(s_{\dim A'}^{B_\Gamma})
\end{split}\]
and similarly for the support of the target. The equation on the support of a
term follows again immediately by induction from the definitions of
\(R^{\Term}\), \(\Var\) and the support. This concludes the definition of
the morphism \(R\). Functoriality of \(R\) and naturality of \(R^{\Type}\) and
\(R^{\Term}\) can be shown by a further mutual induction, using the universal
property of the pushout defining \(R\).

\begin{figure}
  \centering
  \fbox{
    \begin{minipage}{.95\linewidth}
      \begin{mathpar}
        R(\emptyset) = \emptyset \and R_{\Gamma, \emptyset}(\sub{}) = \sub{}
        \and R_{\Delta,(\Gamma,x: A)} (\sub{\gamma,t} ) =
        \sub{R_{\Delta,\Gamma}(\gamma),R^{\Term}_{\Delta,A[\gamma]}(t)} \\
        \begin{tikzcd}[ampersand replacement=\&, column sep = large]
          \Free\sphere^{\dim A} \& \Free\disk^{\dim A + 1} \\
          R(\Gamma) \& R(\Gamma, x:A)
          \ar[from=1-1, to=1-2, hook, "\Free\iota_{\dim A}"]
          \ar[from=2-1, to=2-2, dashed, "\pi_{\Gamma, A}"']
          \ar[from=1-1, to=2-1, "R^{\Type}_\Gamma(A)"']
          \ar[from=1-2, to=2-2, dashed, "v_{\Gamma,A}"]
          \ar[from=1-1, to=2-2, phantom, "\ulcorner"{very near end}]
        \end{tikzcd} \and
        \begin{aligned}
          R_\Gamma^{\Type}(\obj)
          &= \sub{} \\
          R_{\Gamma}^{\Type}(\arr[A]{u}{v})
          &= \sub{R_{\Gamma,A}^{\Term}(u), R_{\Gamma,A}^{\Term}(v)}
        \end{aligned} \\
        R^{\Term}_{(\Gamma, x:A), B}(\var v) =
        \begin{cases}
          v_{\Gamma,A} & \text{if } v=x \\
          \pi_{\Gamma, A} \circ R^{\Term}_{\Gamma,B}(\var v) &\text{otherwise}
        \end{cases} \and
        R^{\Term}_{\Delta,A}(\coh_{\Gamma,A}[\gamma])
        = \coh(B_{\Gamma}, \Free R_\Gamma^{\Bat}\circ R_\Gamma^{\Type}A ,
        R_{\Gamma,\Delta}(\gamma)\circ (\Free R_\Gamma^{\Bat})^{-1})
      \end{mathpar}
    \end{minipage}
  }
  \caption{The morphism \(R\colon \syn_{\catt}\to \Comp^{\op}\)}
  \label{fig:comparison-functor-catt}
\end{figure}

\begin{theorem}\label{thm:main}
  The functor \(R \colon \syn_{\catt}\to \Comp^{\op}\) is fully faithful with
  essential image the finite computads. The natural transformations
  \(R^{\Type}\) and \(R^{\Term}\) are invertible.
\end{theorem}
\begin{proof}
  We first prove faithfulness together with injectivity of the natural
  transformations \(R^{\Type}\) and \(R^{\Term}\) by mutual induction on the
  syntax.
  Suppose that \(R^{\Term}(u) = R^{\Term}(v)\) for terms \(\Gamma\vdash u: A\)
  and \(\Gamma\vdash v:A\). Then if \(u\) is a variable, so is \(v\) and in that
  case \(u=v\). If \(u = \coh_{\Gamma,A}[\gamma]\) is a coherence term, then so
  is \(v = \coh_{\Gamma',A'}[\delta]\) and the above equality gives:
  \begin{align*}
    B_{\Gamma} &= B_{\Gamma'}\\
    \Free (R^{\Bat}_{\Gamma})\circ R^{\Type}(A)
    &= \Free (R^{\Bat}_{\Gamma'})\circ R^{\Type}(A')\\
    R\gamma\circ\Free (R^{\Bat}_{\Gamma})^{-1}
    &= R\gamma'\circ\Free (R^{\Bat}_{\Gamma'})^{-1}.
  \end{align*}
  Using the injectivity of the map \(B_{\bullet}\) and the inductive hypothesis
  for \(R^{\Type}(A)\) and \(R\gamma\), we conclude that \(u=v\). Similarly, one
  can prove that if \(R^{\Type}(A) = R^{\Type}(A')\), then \(A = A'\) by case
  analysis on \(A\). Finally, one can prove that if \(R\gamma = R\delta\), then
  \(\gamma=\delta\) by case analysis on \(\gamma\) and the uniqueness part of
  the universal property of the pushout defining \(R\) on contexts.

  We then show fullness together with surjectivity of \(R^{\Type}\) and
  \(R^{\Term}\) again by mutual induction on dimension, and structural induction
  on the cells and morphisms of computads. Fix a context \(\Delta\vdash\). First
  note that \(R^{\Term}\) is surjective on generator cells since it induces a
  bijection between the variables of \(\Delta\) and the generators of
  \(R\Delta\). Consider an \(n\)\-cell of the form \(c = \coh(B,A,f)\). By
  Theorem~\ref{thm:psctx-trees}, there exists a pasting context \(\Gamma\) such
  that \(B = B_{\Gamma}\). By surjectivity of \(R^{\Type}\) on
  \((n-1)\)-spheres, and by structural induction, there exist a type
  \(\Gamma\vdash A'\), and a substitution \(\Delta\vdash\gamma : \Gamma\)
  respectively, satisfying
  \begin{align*}
    (\Free R^{\Bat}_{\Gamma})^{-1}\circ  A &= R^{\Type}(A') &
    f \circ (\Free R^{\Bat}_{\Gamma}) &= R\gamma .
  \end{align*}
  By the equality of support proven in the construction of \(R\), the type
  \(A'\) must satisfy the side-condition of the rule~\nameref{rule:coh}. By
  construction, \(c = R^{\Term}(\coh_{\Gamma,A'}[\gamma])\), hence \(R^{\Term}\)
  is surjective on \(n\)\-cells. Surjectivity of \(R^{\Type}\) onto
  \(n\)\-spheres follows by surjectivity of \(R^{\Term}\) on \(n\)\-cells.
  Finally, consider a context \(\Gamma\), and a morphism
  \(f : R\Gamma \to R \Delta\). If \(\Gamma\) is the empty context, then \(f\)
  is the unique morphism out of the initial object, which is the image of the
  empty substitution. If \(\Gamma = (\Gamma',x:A)\) is an extension, then \(f\)
  decomposes as \(f = \sub{g,c}\) with \(g : R\Gamma' \to R\Delta\) and \(c\) a
  cell of \(R\Gamma\). By structural induction, \( g = R\gamma\) and
  \(c = R^{\Term}(u)\) for some substitution \(\Delta\vdash \gamma:\Gamma'\) and
  some term \(u\) in \(\Gamma\). By definition of \(R\) on morphisms,
  \(f = R\sub{\gamma,u}\).

  Finally, we show that the essential image of \(R\) coincides with the
  subcategory of finite computads. By construction, the computad \(R\Gamma\) is
  finite with exactly as many generators as the length of \(\Gamma\). We prove
  the converse by induction on the number of generators. If \(C\) is a computad
  with no generators, then \(C = \emptyset = R\emptycontext\). Otherwise, we
  choose a generator \(v \in V_{n}^{C}\) of maximal dimension and consider the
  computad \(C'\) obtained by removing that generator from \(C\). By
  construction, there exists a pushout square
  \[
    \begin{tikzcd}[column sep = large]
      \Free \sphere^{n-1}
      \ar[r,"\Free\iota_{n}"]
      \ar[d,"\phi_{n}^{C}(v)"']
      \ar[rd,phantom,"\ulcorner" very near end]
      & \Free\disk^{n}\ar[d,"v"] \\
      C'\ar[r,"\iota"'] & C
    \end{tikzcd}
  \]
  By the inductive hypothesis, there exists some context \(\Gamma\) together
  with an isomorphism \(f \colon C' \xrightarrow{\sim} R\Gamma\). Then by
  uniqueness of the pushout, we have that
  \[
    C \cong R(\Gamma,x:A),
  \]
  where \(A = (R^{\Type}_{\Gamma})^{-1}(f\circ \phi_{n}^{C}(v))\). This
  concludes the proof that \(R\) is essentially surjective onto finite
  computads.
\end{proof}

\begin{corollary}
  The syntactic category \(\syn_{\catt}\) is equivalent to the opposite of the
  full subcategory of \(\Comp\) consisting of finite computads.
\end{corollary}

\section{Models of CaTT}
\label{sec:models}

We conclude our study with a a brief discussion about the models of the
dependent type theory \catt. Denote \(\PsCtx\) the full subcategory of
\(\syn_{\catt}\) whose objects are the contexts satisfying \(\Gamma\vdashps\).
Explicitly, \(\PsCtx\) is the category whose objects are pasting contexts, and
morphisms are the substitutions in \catt between them, viewing them as contexts
in \catt. Benjamin, Finster and Mimram~\cite{benjamin_globular_2021} have showed
that the models of the dependent theory \catt are equivalent to the copresheaves
over \(\PsCtx\) that preserve a class of limits, called globular products.

On the other hand, the second author's thesis~\cite{markakis_thesis_2024} shows
that the free \(\omega\)\-category monad \(T\) has arities the class of
globular pasting diagrams, which means that the nerve functor
\begin{gather*}
  N_T \colon \operatorname{Alg}_T \to [\Bat^{\op}, \Set] \\
  N_T(X)(B) = \operatorname{Alg}_T(F\Pos B, X)
\end{gather*}
is fully faithful, where \(\Bat\) is the category of Batanin trees
and morphisms of computads, and \(F\) is the free \(T\)\-algebra functor. It is
essentially surjective onto presheaves that preserve a classs of colimits,
called globular sums. Alternatively,
this follows from the comparison theorem of Dean et
al.~\cite[Corollary~7.36]{dean_computads_2022}, and the work of
Berger~\cite[Theorem 1.17]{berger_cellular_2002}.

The functor \(R\) we have defined in Section~\ref{sec:main} satisfies the
following equality for every pasting context \(\Gamma\),
\begin{align*}
  R\Gamma &= \Free V\Gamma = \Free \Pos{B_{\Gamma}}.
\end{align*}
Since \(B_\bullet\) is bijective and the functor \(R\) is fully faithful, the
latter restricts to an isomorphism of categories \(\PsCtx \cong \Bat^{\op}\).
One can check that under this isomorphism globular products correspond to
globular sums. Composing those equivalences, we conclude that the categories
of \(T\)\-algebras and models of \catt are equivalent.

\bibliographystyle{./entics}
\bibliography{bibliography}

\end{document}